\documentclass[10pt]{amsart}
\usepackage{amssymb,amscd}
\numberwithin{equation}{section}

\newtheorem{theorem}{Theorem}[section]

\newtheorem{proposition}[theorem]{Proposition}
\newtheorem{lemma}[theorem]{Lemma}
\theoremstyle{definition}

\theoremstyle{remark}

\renewcommand{\L}{\mathcal{L}}

\newcommand{\U}{\mathcal{U}}

\newcommand{\R}{\mathbb{R}}
\newcommand{\C}{\mathbb{C}}
\newcommand{\Z}{\mathbb{Z}}

\newcommand{\bbT}{\mathbb{T}}
\newcommand\lie[1]{\mathfrak{#1}}

\newcommand{\fg}{\lie{g}}

\newcommand{\ft}{\lie{t}}

\def	\inv	{^{-1}}

\newcommand\coker{\mathop{\rm coker}\nolimits}
\newcommand\codim{\mathop{\rm codim}\nolimits}

\newcommand\vo{\mathaccent23V}

\begin{document}

\title{ Homotopy  groups of $K$-contact toric manifolds}

\author{Eugene Lerman}
\address{Department of
Mathematics, University of Illinois, Urbana, IL 61801}
\email{lerman@math.uiuc.edu}

\thanks{Supported by the NSF grant DMS-980305.}
\date{\today}

\begin{abstract} 
We compute the first and second homotopy groups of a class of contact
toric manifolds in terms of the images of the associated moment map.
\end{abstract}

\maketitle

\section{Introduction}

In this paper I compute the first and second homotopy groups of
certain toric symplectic cones or, equivalently, of certain contact toric
manifolds.  The main result of the paper is Theorem~\ref{thm1} (the
terms used in the statement are explained below):

\begin{theorem}\label{thm1}
Let $G$ be a torus with Lie algebra $\fg$ and integral lattice $\Z_G =
\ker \{ \exp: \fg \to G\}$.  Let $(B, \xi = \ker \alpha)$ be a contact toric
$G$-manifold of Reeb type with moment cone $C\subset \fg^*$, which is
a strictly convex rational polyhedral cone.  Let $\L$ denote the
sublattice of $\Z_G$ generated by the normal vectors to the
facets of $C$.  The fundamental group of $B$ is the finite abelian
group $\Z_G/\L$.  The second homotopy group of $B$ is a free abelian
group of rank $N - \dim G$ where $N$ is the number of facets of the
cone $C$.
\end{theorem}

Let us recall the necessary definitions (see \cite{CTM} for more
details; see also \cite{LS}).  A manifold $B$ with a contact structure
$\xi = \ker \alpha$ ($\alpha$ is a contact form) is a {\bf toric
$G$-manifold} if there exists an effective action of a torus $G$ on
$B$ preserving $\xi$ with $\dim B + 1 = 2 \dim G$.  By averaging over
the group, if necessary, one can always assume that the torus $G$
preserves a contact form $\alpha$ defining $\xi$.
Given an action of a Lie group $G$ on a manifold $B$ preserving a
contact form $\alpha$, the corresponding {\bf $\alpha$-moment map}
$\Psi_\alpha :B \to \fg^*$ ($\fg^*$ denotes the vector space dual of
the Lie algebra $\fg$ of $G$) is defined by
$$
\langle \Psi_\alpha (b) , X \rangle = \alpha _b (X_B (b))
$$ 
for all $b\in B$, all $X\in \fg$.  As usual $\langle \cdot , \cdot
\rangle$ denotes the canonical paring between $\fg^*$ and $\fg$, and
$X_B$ denotes the vector field on $B$ induced by $X$.  If $f \in
C^\infty (B)^G$ is an invariant function, then $\alpha' = e^f \alpha$
is another contact form defining the same contact distribution $\xi
$ as $\alpha$.  Clearly $\Psi_{e^f \alpha} = e^f \Psi_\alpha$, so the
moment map is an invariant of the contact form and not of the contact
distribution.  On the other hand the subset $C(\Psi)= C(\Psi_\alpha)$
of $\fg^*$ for an $\alpha$-moment map $\Psi_\alpha :B \to \fg^*$
defined by 
$$ 
C(\Psi) = \{ t\Psi_\alpha (b) \mid t\geq 0, \, b \in B\}
$$ 
depends only on the action of $G$ on $B$ and on the contact
distribution $\xi$ but not on the contact form $\alpha$ {\em
per se}.  We will refer to $C(\Psi)$ as the {\bf moment cone} of
the action.
Since a moment map $\Psi_\alpha :B\to \fg^*$ completely encodes the
action of $G$ on $(B, \alpha)$ we regard a contact toric $G$-manifold
as a triple $(B, \xi = \ker \alpha, \Psi_\alpha: B\to \fg^*)$.  Note
that the {\bf symplectization} $(M, \omega) := (B\times \R, d(e^t
\alpha))$ $(t\in \R)$ is a noncompact symplectic toric manifold with moment 
map $\Phi (b,t) = e^t \Psi_\alpha (b)$.  Remark that $\Phi (M) \cup
\{0\} = C(\Psi)$.  Conversely, if a symplectic toric $G$-manifold $(M,
\omega, \Phi :M\to \fg^*)$ is a {\bf symplectic cone}, i.e., 
if there is a free proper action $\{\rho_t\}$ of $\R$ on $M$ commuting
with the action of $G$ such that $\rho_t^* \omega = e^t \omega$, then
$M/\R$ is naturally a contact toric manifold.

A contact manifold $(B, \xi = \ker \alpha)$ with an action of a torus
$G$ preserving $\alpha$ is {\bf Reeb type} if there is $X\in \fg$ such
that the function $\langle \Psi_\alpha, X\rangle = \iota (X_B) \alpha$
is strictly positive.
By a result of Boyer and Galicki \cite{BG} (see also Theorem~4.3 in
\cite{LS}), the moment cone of a contact toric $G$ manifold of Reeb
type is a strictly convex rational polyhedral cone.  ``Strictly
convex'' means that the moment cone contains no linear subspaces of
positive dimension, i.e., it's a cone on a polytope.  ``Rational
polyhedral'' means that there exist vectors $\mu_1, \ldots,
\mu_N$ in the integral lattice $\Z_G := \ker (\exp :\fg \to G)$ of the
torus $G$ such that 
$$ 
C(\Psi) = 
\{\eta \in \fg^* \mid \langle \eta, \mu_j\rangle \geq 0, \, j=1, \ldots, N\}.  
$$

There are several reasons for wanting to compute the homotopy groups of contact
toric manifolds of Reeb type.  \\[5pt]
\noindent
1. All contact manifolds of Reeb type are $K$-contact (see
Proposition~\ref{reeb} below), hence the title of the paper.  In fact
contact {\em toric} manifolds of Reeb type are Sasakian, as proved by
Boyer and Galicki (Theorem~5.3 in \cite{BG}).  Methods recently
developed by Boyer, Galicki, Mann and others use Sasakian structures
to obtain explicit positive Einstein metrics. \\[5pt]
\noindent
2. A classification of contact toric manifolds \cite{CTM} shows that
contact toric manifolds not diffeomorphic to the ones of Reeb type are
easy to understand: they are either $S^2 \times S^1$, or products
$\bbT^k \times S^{k+2l -1}$ ($k>1$, $l\geq 0$) or principal $\bbT^3$
bundles over $S^2$.  So if one wants to understand the {\em topology}
of contact toric manifolds, the manifolds of Reeb type are the ones to
concentrate on. \\[5pt]
\noindent
3. One motivation for studying the topology of contact toric manifolds
is their apparent difference from (topological) toric
manifolds. Recall that in 1991 Davis and Januszkiewiecz defined
(topological) toric manifolds as manifolds with torus action locally
modeled on the standard action of $\bbT^n$ on $\C^n$ and having a
simple polytope as the orbit space \cite{DJ}.  Such a manifold is
determined by a polytope and a characteristic function, a function
that assigns a 1-parameter subgroup of the torus to every facet of
the polytope. They proved a beautiful formula for the integral
cohomology ring of a toric manifold; it is the Stanley-Reisner ring of
the polytope modulo an ideal determined by the characteristic function
(for smooth projective toric varieties the formula is known as the
Danilov-Jurkiewicz theorem).  In particular the cohomology ring is
generated by elements of degree two, odd dimensional cohomology
vanishes and there is no torsion.  They also proved that such
manifolds are simply connected.  In contrast, the odd dimension
cohomology of a contact toric manifold need not vanish (cf.\ $\R
P^3$), there is torsion and the fundamental group need not be trivial.
\\[5pt]
\noindent
4. Another motivation comes from the study of completely integrable
geodesic flows.  According to Toth and Zelditch \cite{TZ}, a geodesic
flow on a manifold $Q$ is {\bf toric integrable} if there exists a
homogeneous completely integrable action of a torus on the punctured
cotangent bundle $T^*Q \smallsetminus Q$ which preserves the geodesic
flow.  Naturally in this case the co-sphere bundle $S(T^*Q)$ is a
contact toric manifold.  It would be interesting to find a topological
obstruction to the existence of a toric integrable geodesic flow on a
compact manifold $Q$ and for that one needs to understand the topology
of contact toric manifolds.\\

\noindent
We now outline the proof of  Theorem~\ref{thm1}.\\

1) Since a contact manifold $B$ is homotopy equivalent to its
symplectization $M = B\times \R$, we compute the homotopy groups of
the symplectization.\\

2) The symplectization $M$ of $B$ is the symplectic quotient at 0 of $\C^N
\smallsetminus \{0\}$ by a compact abelian group $T$ with 
$\pi_0(T) = \Z_G /\L$ and $\dim T = N -\dim G$.  That is to say, $M =
(\Phi_T\inv (0)\smallsetminus \{0\})/T$ where $\Phi_T : \C^N \to
\ft^*$ denotes the $T$-moment map $\Phi_T : \C^N \to \ft^*$ for the linear action of $T$ on $\C ^N$.\\

3) The set $\Phi_T\inv (0)\smallsetminus \{0\}$ has the homotopy type
of $\C^N \smallsetminus (V_1 \cup V_2 \cup \ldots \cup V_r)$ where
each $V_j \subset \C^N$ is a linear subspace of complex codimension at
least 2.  Hence $\pi_0 (\Phi_T\inv (0)\smallsetminus \{0\})=\pi_1
(\Phi_T\inv (0)\smallsetminus \{0\})=\pi_2 (\Phi_T\inv
(0)\smallsetminus \{0\})= *$.\\

4) Since the group $T$ acts freely on $\Phi_T\inv (0)\smallsetminus
\{0\}$, we see from the long exact sequence of homotopy groups for the
fibration $T\to (\Phi_T\inv (0)\smallsetminus \{0\}) \to M$ that
$$
\pi_1 (M) = \pi_0 (T) \quad \text{and}\quad \pi_2 (M) = \pi_1 (T).
$$ 
The details of the argument are the subject of the next section.  In
the last section we explain the connection between torus actions of
Reeb type and being $K$-contact.

\subsection*{A note on notation}  

Throughout the paper the Lie algebra of a Lie group denoted by a
capital Roman letter will be denoted by the same small letter in the
fraktur font: thus $\fg$ denotes the Lie algebra of a Lie group $G$
etc.  The natural pairing between $\fg$ and its vector space dual
$\fg^*$ is denoted by $\langle \cdot, \cdot \rangle$.  If $A: V \to W$
is a linear map, we denote the corresponding map on the dual spaces by
$A^*$, $A^*: W^* \to V^*$.

When a Lie group $G$ acts on a manifold $M$ we denote the action by an
element $g\in G$ on a point $x\in M$ by $g\cdot x$; $G\cdot x$ denotes
the $G$-orbit of $x$ and so on.  The vector field induced on $M$ by an
element $X$ of the Lie algebra $\fg$ of $G$ is denoted by $X_M$.  Thus
$X_M (m) =\left. \frac{d}{dt}\right| _{t=0} (\exp tX)\cdot m$.

For us a {\em torus} is a compact connected abelian group.  If $G$ is
a torus, we denote its weight lattice by $\Z_G^*$, it is a subgroup of
$\fg^*$.  The dual lattice of $\Z_G^*$ is the integral lattice $\Z_G$.
Recall that $\Z_G = \ker (\exp :\fg \to G)$.  Thus $G = \fg/\Z_G$.

\subsection*{Acknowledgments} I thank Charles Boyer, Sue Tolman and Bill 
Graham for a  number of useful conversations. 

\section{Proof of the main result, \protect Theorem~\ref{thm1}}

It was proved in \cite{CTM} that the moment cone $C(\Psi)$ of a
(compact connected) contact toric $G$-manifold $(B, \xi = \ker \alpha,
\Psi_\alpha: B\to \fg^*)$ of Reeb type is a {\em good cone}.  This means the 
following.  Let $\{F_i\}$ denote the set of facets (codimension one
faces) of $C(\Psi)$.  Since $C(\Psi)$ is rational, each facet is of
the form 
$$ 
F_i = \{ \eta \in C(\Psi) \mid \langle \eta, \mu_i\rangle= 0\} 
$$ 
for some primitive vector $\mu_i$ in the integral lattice
$\Z_G$ of $G$.  Then
\begin{enumerate}
\item every codimension $\ell$, $0<\ell < \dim G$, face $F$ of $C(\Psi)$ can 
be written uniquely as
$$
F = F_{i_1} \cap \ldots F_{i_\ell}
$$
where $F_{i_j}$'s are the facets containing $F$, and 
\item  the $\Z$-module generated by the normals to the facets 
$F_{i_1}, \ldots, F_{i_\ell}$ is a direct summand of $\Z_G$ of rank $\ell$.
\end{enumerate}

We have a uniqueness result \cite{CTM}: if $(B, \xi = \ker \alpha,
\Psi_\alpha)$ and $(B', \xi'=\ker \alpha', \Psi_{\alpha'})$ are two
(compact connected) contact toric manifolds of Reeb type and the
moment cones are equal then the contact toric manifolds are
equivariantly contactomorphic.

There is also a corresponding existence result. Given a good polyhedral
cone $C\subset \fg^*$ (where $\fg^*$ is the dual of the Lie algebra of
a torus $G$) there exists a compact connected contact toric $G$-manifold
$(B_C, \xi_C = \ker \alpha_C,
\Psi_{\alpha_C})$ with the moment cone $C(\Psi_C)$ equal to $C$
(Theorem~2.18(4) of \cite{CTM}).  Moreover $(B_C, \xi_C = \ker
\alpha_C, \Psi_{\alpha_C})$ can be constructed as a contact quotient of 
the standard odd dimensional sphere.  In fact it is more convenient to
construct the symplectization $(M_C, \omega_C, \Phi_C: M_C \to \fg^*
)$ of $(B_C, \alpha_C, \Psi_{\alpha_C} : B_C \to \fg^*)$.  Then for
any contact toric $G$-manifold $(B', \xi'=\ker \alpha',
\Psi_{\alpha'})$ with $C(\Psi_{\alpha'}) =C$ we have 
$$
\pi_1 (M_C) = \pi_1 (B'), \quad \pi_2 (M_C) = \pi_2 (B')
$$ and so on.  Note that the moment map image $\Phi_C (M_C)$ is
$C\smallsetminus \{0\}$.

Recall from \cite{CTM} the construction of the symplectic toric
manifold $(M_C, \omega_C, \Phi_C:M_C \to \fg^*)$.  As above let
$\mu_1, \ldots, \mu_N \in \Z_ G$ denote the primitive inward normals
to the facets of the good strictly convex cone $C$.  Since $C$ is
strictly convex and has non-empty interior, span$_\R \{\mu_i\} =\fg$.
Hence the abelian group $\Z_G/\L$, where $\L = \text{span}_\Z
\{\mu_i\}$, is finite.  Consider the $\Z$-linear map 
\begin{equation}
\varpi: \Z^N \to \Z_G, \quad
\varpi (a_1, \ldots, a_N) = \sum a_i \mu_i. 
\end{equation} 
Its cokernel is $\Z_G/\L$.  It extends to a surjective $\R$-linear map
\begin{equation} \label{eq2}
\tilde{\varpi}: \R^N \to \fg, \quad 
\tilde{\varpi} (a_1, \ldots, a_N) = \sum a_i \mu_i,
\end{equation}
which drops down to a surjective Lie
group homomorphism 
$$
\bar{\varpi} : \bbT^N = \R^N/\Z^N \to \fg/\Z_G =
G,  \quad 
\bar{\varpi} ([a_1, \ldots, a_N]) = \exp (\tilde{\varpi} (a_1,
\ldots, a_N)) = \exp (\sum a_i \mu_i).
$$  
Here $[a_1, \ldots, a_N]$ denotes the class of $(a_1, \ldots, a_N) \in
\R^N$ in $\bbT^N$ and $\exp : \fg \to G$ denotes the exponential map.
Let $T = \ker \bar{\varpi}$; it is a closed by not necessarily
connected subgroup of $\bbT^N$.  The standard linear action of
$\bbT^N$ on $\C^N$ preserving the standard symplectic form $ \sqrt{-1}
\sum dz_j \wedge d\bar{z}_j$ gives rise to a linear symplectic action
of $T \subset \bbT^N$.  Denote the corresponding homogeneous moment
map by $\Phi_T$; $\Phi_T: \C^N \to
\ft^*$.  The moment map $\Phi: \C^N \to (\R^N)^*$ for the standard action of 
$\bbT^N$ on $\C^N$ is given by the formula 
\begin{equation}\label{eq-star}
\Phi (z_1, \ldots, z_N) = \sum |z_j|^2 e_j^*
\end{equation}
where $e_1^*, \ldots, e_N^*$ is the standard basis
of $(\R^N)^*$.  Hence, if $\iota : \ft \to \R^N$ denotes the inclusion
of the Lie algebra $\ft$ of $T$, we have
$
\Phi_T = \iota^* \circ \Phi.
$
We recall from \cite{CTM}:
\begin{lemma} \label{lem0}
We use the notation above.  The set $\Phi_T \inv (0)\smallsetminus
\{0\}$ is a manifold.  The group $T$ acts freely on this manifold. The
symplectic manifold $M:= (\Phi_T \inv (0)\smallsetminus \{0\})/T$ is
the desired $G= \bbT^N/T$ symplectic manifold, that is, it is a
symplectic cone and the image of the $G$-moment map is $C\smallsetminus
\{0\}$. In particular $\Phi (\Phi_T \inv (0)) = \tilde{\varpi}^* (C)$ where
$\tilde{\varpi}^*: \fg^* \to (\R^N)^*$ is dual to $\tilde{\varpi}$ 
(cf.\ (\ref{eq2})).
\end{lemma}

Our proof of Theorem~~\ref{thm1} is based on two lemmas.  The first
one describes the group $\pi_0 (T)$ of connected components of $T$:

\begin{lemma} \label{lem1}
Let $T \subset \bbT^N$ be as above.  Then $\pi_0 (T) = \Z_G/\L$ where,
as above, $\L$ is the sublattice of the integral lattice $\Z_G$
spanned by the primitive normals to the facets of the cone $C$.
\end{lemma}
The second lemma shows that the manifold 
$\Phi_T \inv (0)\smallsetminus \{0\}$ has the homotopy type of $\C^N
\smallsetminus (V_1 \cup \ldots \cup V_r)$ where $V_j \subset \C^N$
are complex linear subspaces of complex codimension at least 2.  In
fact the subspaces $V_j$ being deleted are determined by the
combinatorics of the polyhedral cone $C$.  To make this precise we
need a few definitions.
For a subset $I\subset \{1, \ldots, N\}$ define the corresponding
coordinate subspace $V_I$ by 
$$ 
V_I : = \{ z\in \C^N \mid j\in I\Rightarrow z_j = 0\} = 
\bigcap _{j\in I} \{z_j = 0\}.  
$$ 
For each $j\in \{1, \ldots, N\}$ the $j$th facet $F_j$ of the cone $C$
satisfies
$$
F_j = C\cap \{\eta\in \fg^* \mid \langle \eta, \mu_j \rangle = 0\}.
$$
Now consider the set 
$$
\U := \left\{ I \subset \{1, \ldots, N\} \mid \bigcap _{j\in I} F_j = \{0\} 
\right\},
$$ 
the collection of subsets $I$ of $ \{1, \ldots, N\}$ such that the
facets indexed by the elements of $I$ intersect only at the vertex.
\begin{lemma} \label{lem2}  
The manifold $\Phi_T \inv (0)\smallsetminus \{0\}$ has the same 
homotopy type as
\begin{equation}
\C ^N \smallsetminus \bigcup _{I\in \U} V_I.
\end{equation}
\end{lemma}

Let us assume the lemmas for a moment and prove the main theorem,
Theorem~\ref{thm1}.  
\begin{proof}[Proof of Theorem~\ref{thm1}]
As was remarked previously, it is enough to prove that the symplectic
toric manifold $M_C = M = (\Phi_T\inv (0)\smallsetminus \{0\})/T$ has
the properties that $\pi_1 (M) = \Z_G /\L$ and that $\pi_2 (M_C) =
\Z^d$ where $d = N - \dim G$.  Since $T$ acts freely on 
$Z:= \Phi_T\inv (0)\smallsetminus \{0\}$, we have a long exact
sequence of homotopy groups
\begin{equation}\label{eq-long-exact}
\cdots \to \pi_2 (Z) \to \pi_2 (M) \to \pi _1 (T) \to \pi_1(Z) 
\to \pi_1 (M) \to \pi_0 (T) \to  \pi_0 (Z) \to \pi_0 (M).
\end{equation}
Since every facet $F_j$ of $C$ is not $\{0\}$, the set $\U$ contains
no singletons.  Since $\dim _\C V_I = N - |I|$, it follows that for
any $I\in \U$, $\codim_\C V_I = |I| \geq 2$. Hence by Lemma~\ref{lem2}
$Z$ is connected and the homotopy  groups $\pi_2(Z)$, $\pi_1 (Z)$ are trivial.
It follows from (\ref{eq-long-exact}) that 
$$
\pi_2 (M) = \pi_1 (T) \text{ and } \pi_1 (M) = \pi_0 (T).
$$  
By Lemma~\ref{lem1} $\pi_0 (T) = \Z_G/ \L$.  Clearly $\pi_1 (T) = \Z^d$,
$d= \dim T = \dim \bbT^N -\dim G$.  
\end{proof}
\begin{proof}[Proof of Lemma~\ref{lem1}]
This is a simple application of Snake lemma.  Consider the commuting diagram
\begin{equation*}
\begin{CD}
0 @>>> \Z^N @>>> \R^N @>\exp >> \bbT^N @>>> 1\\
@.  @VV{\varpi}V @VV{\tilde{\varpi}}V @VV{\bar{\varpi}}V\\
0 @>>> \Z_G @>>>\fg @>\exp >> G @>>> 1
\end{CD}
\end{equation*}
By Snake lemma we have a long exact sequence 
$$
\ker \varpi \to \ker \tilde{\varpi} \to \ker \bar{\varpi} \to \coker \varpi
 \to \coker \tilde{\varpi} \to \coker \bar{\varpi} .  
$$ 
By construction $\tilde{\varpi}$ is onto, hence $\coker \tilde{\varpi}
=0$.  On the other hand $\coker \varpi = \Z_G /\L$.  By definition
$\ker \bar{\varpi} = T$, $\ker \tilde{\varpi} =\ft$ and the map $\ker
\tilde{\varpi} \to \ker \bar {\varpi}$ is simply the exponential map
$\exp: \ft \to T$.  Since $\coker (\exp :\ft \to T)$ is $\pi_0 (T)$ we
get $\pi_0 (T) \simeq \Z_G/\L$.
\end{proof}

\begin{proof}[Proof of Lemma~\ref{lem2}]
We keep the notation of the discussion above.  The proof is an
elementary application of the correspondence between symplectic
quotients and Geometric Invariant Theory (GIT) quotients as developed
by Mumford, Guillemin, Sternberg, Kirwan, Neeman, Sjamaar and others.
The key point is that the GIT quotient $\C ^N /\!/T^\C $ and the
symplectic quotient $\Phi_T\inv (0) /T$ are isomorphic as stratified
spaces.  It will be most convenient for us to quote \cite{Sj} where
Kirwan's results on the isomorphism between symplectic and GIT
quotients were suitably refined.

\noindent
(1)\ By Lemma~\ref{lem0} the group $T$ acts freely on the manifold $Z
= \Phi_T \inv (0) \smallsetminus \{0\}$.\\ 
\noindent
(2) \ By Example~2.3 of \cite{Sj}, $\Phi_T$ is admissible in the sense
of \cite{Sj} p.~109, and the set of analyticly semi-stable points
$(\C^N)^{ss}$ for the action of $T$ on $\C^N$ is all of $\C^N$.\\
\noindent
(3) \ By Proposition~1.6 of \cite{Sj} for any point $z\in \C^N$ the
stabilizer in the complexified group is the complexification of the
stabilizer: 
$$
(T^\C)_z = (T_z)^\C .
$$
Hence by (1), $(T^\C)_z$ is trivial for all $z\in Z$.\\
\noindent
(4) \ By Proposition~2.4(ii) of \cite{Sj} the orbit $T^\C \cdot z$ is
closed in $(\C ^N)^{ss} = \C ^N$ if and only if $T^\C \cdot z \cap
\Phi_T\inv (0) \not = \emptyset$.  Thus
\begin{equation} \label{1*}
\{z\in \C ^N \mid T^\C \cdot z \text{ is closed }\} = 
	\{0\} \cup T^\C \cdot Z .
\end{equation}
\noindent
(5) \ Since the actions of $(\bbT^N)^\C$ and $T^\C$ commute, the union
(\ref{1*}) of closed $T^C$ orbits is $(\bbT^N)^\C$ invariant.  Hence, since 
$\{0\}$ is fixed by $(\bbT^N)^\C$, the set 
$$
S := T^\C \cdot Z
$$
is $(\bbT^N)^\C$ invariant.\\
\noindent
(6) \ Proposition~2.4(iii) of \cite{Sj} implies that $(T^\C \cdot
Z)/T^\C = Z/T$.  Combining this with (3) we see that $S$ is a
$T^\C/T$-bundle over $Z$.  Since $T^\C/T$ is diffeomorphic to the Lie
algebra $\ft$ of $T$, the manifolds $S$ and $Z$ are homotopy
equivalent.\\
\noindent
(7) \ For any subset $I$ of $\{1, \ldots, N\}$ define
$$
\vo_I = \{ z\in V_I \mid z_j \not = 0 \text{ for } j\not \in I\},
$$ the ``interior'' of the coordinate subspace $V_I$.  The set $\vo_I$
is a single $(\bbT^N)^\C$ orbit.  It satisfies
$$
\vo_I = V_I \smallsetminus \bigcup _{I' \supset I, I' \neq I}V_{I'}. 
$$
We claim that
\begin{equation}\label{2*}
\vo_I \subset S \Leftrightarrow 
	\bigcap _{j\in I} F_j \text{ is a nonzero face of } C.
\end{equation} 
\begin{proof}[Proof of (\ref{2*})]
Note that since $S$ is $(\bbT^N)^\C$ invariant and $\vo_I$ is a
$(\bbT^N)^\C$ orbit, $\vo_I \subset S$ $\Leftrightarrow$ $\vo_I \cap S
\not = \emptyset$.  Also, since $z\in S \Leftrightarrow T^\C \cdot z
\cap Z \not = \emptyset$ and since $S$ is $(\bbT^N)^\C$ invariant, we
have $$ z\in S \Leftrightarrow (\bbT^N)^\C \cdot z \cap Z \not =
\emptyset.  $$ As before let $\mu_j\in \Z_G$ denote the (primitive
inward pointing) normal to the facet $F_j$ of $C$.  Suppose $F_I := \bigcap
_{j\in I} F_j$ is a nonzero face of $C$.  Pick  a point $\eta$ in
the relative interior of $F_I$. Then $\langle \eta, \mu_k \rangle >
0$ for all $k \not \in I$.  Let $z^\eta_j = \sqrt{ \langle \eta,
\mu_j\rangle}$; $z^\eta := (z_1^\eta, \ldots, z_N^\eta)$
satisfies 
$$
\langle \Phi (z^\eta), e_j\rangle = |z^\eta_j|^2 = \langle \eta, \mu_j\rangle
= \langle \eta, \tilde{\varpi} (e_j)\rangle = \langle
\tilde{\varpi}^* (\eta), e_j\rangle 
$$ for all $j$, where, as before, $e_1, \ldots, e_N$ is the standard
basis of $\R^N$, $\Phi: \C^N \to (\R^N)^*$ is the moment map for the
standard action of $\bbT^N$ on $\C^N$ (see (\ref{eq-star})) and
$\tilde{\varpi}: \R^N \to \fg$ is the surjective map defined earlier
by (\ref{eq2}).  Hence $\Phi (z^\eta) = \tilde{\varpi}^* (\eta)$, so
$z^\eta \in \Phi\inv (\tilde{\varpi}^* (eta))$.  Since $\eta \not = 0$
we have 
$$
\emptyset \not =  T^\C \cdot z^\eta \cap (\Phi\inv (\tilde{\varpi}^* (C)) 
\smallsetminus \{0\}) = T^\C \cdot z^\eta \cap Z,
$$ 
where we used the fact that $\Phi \inv (\tilde{\varpi}^* (C)) =
Z\cup \{0\}$. Also $z^\eta \in \vo_I$ since $|z_j^\eta|^2 = \langle \eta,
\mu_j\rangle$ for all $j$ and $\langle \eta, \mu_j\rangle >0$ for 
$j\not \in I$.  This proves that if the intersections $\bigcap _{j\in
I} F_j$ is a nonzero face of $C$ then $\vo_I \cap Z \not = \emptyset
$.  Hence $\vo_I \cap S \not = \emptyset$ and therefore $\vo_I \subset
S$.

Conversely, suppose $\vo_I \subset S$.  Then $\vo_I \cap Z \not =
\emptyset $. For any $z\in \vo_I \cap Z $ we have:
$\Phi (z) \in
\tilde{\varpi}^* (C)$, $|z_j|^2 \not = 0$ for $j\not\in I$, $|z_j|^2 = 0$ 
for $j\in I$.  Therefore $\Phi (z) = \tilde{\varpi}^* (\eta) $ for some
$\eta \in C$ and $\langle \eta, \mu_j\rangle \not = 0$ for all
$j\not\in I$, $\langle \eta, \mu_j\rangle = 0$ for all $j\in I$.
Hence 
$$
\eta\in   \left(\bigcap _{j\not\in I} 
\{ \eta \in \fg^* \mid \langle \eta, \mu_j\rangle  > 0 \}\right)
\cap \left(\bigcap _{j\in I} 
\{ \eta \in \fg^* \mid \langle \eta, \mu_j\rangle  = 0 \}\right).
$$   Thus $F_I = \bigcap _{j\in I} F_j$ is a
nonzero face of $C$.  This proves (\ref{2*}).
\end{proof}
\noindent
(8) \ If $\bigcap _{j\in I} F_j = \{0\}$ then for any $I' \supset I$,
$\bigcap _{j\in I'} F_j = \{0\}$ as well.  Since $V_I = \bigcup
_{I'\supseteq I} \vo_{I'}$, (\ref{2*}) implies that
$$
S = \C^N \smallsetminus \bigcup _{I\in \U} V_I. 
$$ 
By (6) $Z= \Phi_T\inv (0) \smallsetminus \{0\}$ is homotopy equivalent to 
$S=  \C^N \smallsetminus \bigcup _{I\in \U} V_I$ and the result follows.
\end{proof}

\section{Reeb type and $K$-contact}
In this section we prove a version of Proposition~2.1 of Yamazaki
\cite{Yamazaki1} that relates torus actions and $K$-contactness.
Recall that the Reeb vector field $R_\alpha$ on a contact manifold
$(B, \alpha)$ is the unique vector field defined  by the equations
$$
\iota (R_\alpha) d\alpha = 0, \quad \iota (R_\alpha) \alpha = 1.
$$
The Reeb vector field defines a splitting of the tangent bundle of
$B$: 
\begin{equation}\label{eq3.1}
TB = \xi \oplus \R R_\alpha,
\end{equation}
 where $\xi = \ker \alpha$ is the
contact distribution.  Since  $(\xi, d\alpha |_\xi)$ is a symplectic
vector bundle, there exists a complex structure $J$ on $\xi$
compatible with $d\alpha |_\xi$ so that $g_\xi = d\alpha |_\xi (
\cdot, J \cdot)$ is a metric on $\xi$.  Using (\ref{eq3.1}) we may extend 
$g_\xi$ by zero to all of $TB$.  Then $g = g_\xi \oplus \alpha \otimes
\alpha$ is a Riemannian metric on $B$ in which $\xi$ and $R_\alpha$
are orthogonal and the length of the Reeb vector field is 1.  The
metric $g$ is  said to be {\bf adapted} to the contact form
$\alpha$.  If additionally the Reeb vector field is Killing with
respect to an adapted metric $g$, i.e., if $L_{R_\alpha} g = 0$, then
the pair $(\alpha, g)$ is called a {\bf $K$-contact structure} on $B$.
If given a contact distribution $\xi$ on a manifold $B$ there exists a
$K$-contact structure with $\ker \alpha = \xi$ we will say that $(B,
\xi)$ {\bf admits a $K$-contact structure}.

Note that if a Lie group $G$ acts on $B$ preserving a contact form
$\alpha$ then it preserves the Reeb vector field $R_\alpha$, the contact
distribution $\xi = \ker \alpha$ and the symplectic structure $d\alpha
|_\xi$.  Therefore if $G$ is compact we may choose the complex
structure $J$ (and hence the adapted metric $g$) to be $G$-invariant.
 
\begin{proposition}\label{reeb}
 A compact contact manifold $(B,
\xi = \ker \alpha)$ admits the structure of a $K$-contact manifold if and 
only if there exists an action of a torus $G$ on $B$ preserving
$\alpha$ and a vector $X\in \fg$ such that the function $\iota (X_B)
\alpha = \langle \Psi_\alpha , X\rangle$ is strictly positive, i.e., the $G$ 
action is of Reeb type.  Here as before $X_B$ denotes the vector field
on $B$ induced by $X \in \fg$ and $\Psi_\alpha$ denotes the
$\alpha$-moment map.
\end{proposition}

\begin{proof}
Suppose the action of a torus $G$ on $(B, \xi = \ker \alpha)$ is of Reeb type,
i.e., suppose there is a vector $X\in \fg$ such that $\langle
\Psi_\alpha, X \rangle$ is strictly positive (note that this is a condition 
on the co-oriented contact distribution $\xi$ and not just on the
contact form $\alpha$).  We then can multiply $\alpha$ by a positive
$G$-invariant function $f$ so that $\langle \Psi_{f\alpha} , X \rangle
= 1$ ( take $f = 1/\langle\Psi_{\alpha} , X \rangle$).  Therefore it
is no loss of generality to assume that $\alpha (X_B) = \langle
\Psi_\alpha, X \rangle = 1$.  Since $G$ action preserves $\alpha$, we
have $0 = L_{X_B} \alpha = d \iota (X_B) \alpha + \iota (X_B) d \alpha
= d 1 + \iota (X_B) d \alpha$.  Therefore $X_B$ is the Reeb vector
field of $\alpha$.  Now choose an $G$-invariant metric $g$ adapted to
$\alpha$.  Then, since $\alpha$ is $G$-invariant, $L_{X_B} g = 0$, and
so $(\alpha, g)$ is a $K$-contact structure on $(B, \xi)$.

Conversely suppose $(\alpha, g)$ is a $K$-contact structure on $B$.
Since $B$ is compact, the group of isometries of $(B, g)$ is a compact
Lie group $H$.  Take the closure inside $H$ of the flow of the Reeb
vector field $R_\alpha$.  The closure is a compact abelian group $G$,
i.e., a torus.  Since the flow of $R_\alpha$ preserves the contact
form $\alpha$, the action of $G$ preserves $\alpha$ as well.  By
construction $R_\alpha = X_B$ for some vector $X$ in the Lie algebra
of $G$.  Since $R_\alpha$ is a Reeb vector field we have $1 = \iota
(R_\alpha) \alpha = \langle \Psi_\alpha , X \rangle$, where 
$\Psi_\alpha :B\to \fg^*$ is the moment map for the action of $G$ on
$(B, \alpha)$.  Hence the action of $G$ on $(B, \xi = \ker \alpha)$ is
of Reeb type.
\end{proof}

\end{document}